\documentclass[12pt,twoside,reqno]{amsart}
\usepackage{amsmath, amsthm, amscd, amsfonts, amssymb, graphicx, color}
\usepackage[utf8]{inputenc}
\usepackage{cleveref}
\crefname{section}{§}{§§}
\Crefname{section}{§}{§§}
\let \oldsection
\section
\renewcommand{\section}{\vspace{8pt plus 4pt}\oldsection}

\newcommand{\beqa}{\begin{eqnarray*}}
\newcommand{\eeqa}{\end{eqnarray*}}
\newcommand{\beqn}{\begin{eqnarray}}
\newcommand{\eeqn}{\end{eqnarray}}

\newcommand{\R}{\mathbb R}

\newcommand{\mcB}{\mathcal B}

\newcommand{\G}{\Gamma}

\usepackage{lipsum}
\usepackage[super]{nth}
\usepackage{xcolor}
\usepackage{tikz}
\usepackage{pgfplots}
\usepackage{mathrsfs}
\usepackage{graphics}
\usepackage{graphicx} 
\usepackage{epsfig}
\usepackage{wrapfig}
\definecolor{olive}{rgb}{0.3, 0.4, .1}
\definecolor{fore}{RGB}{249,242,215}
\definecolor{back}{RGB}{51,51,51}
\definecolor{title}{RGB}{255,0,90}
\definecolor{dgreen}{rgb}{0.,0.6,0.}
\definecolor{gold}{rgb}{1.,0.84,0.}
\definecolor{JungleGreen}{cmyk}{0.99,0,0.52,0}
\definecolor{BlueGreen}{cmyk}{0.85,0,0.33,0}
\definecolor{RawSienna}{cmyk}{0,0.72,1,0.45}
\definecolor{Magenta}{cmyk}{0,1,0,0}

\newtheorem{thm}{Theorem}[section]

\newtheorem{corollary}[thm]{Corollary}
\newtheorem{lemma}[thm]{Lemma}
\newtheorem{prop}[thm]{Proposition}

\theoremstyle{definition}
\newtheorem{defn}{Definition}[section]

\theoremstyle{remark}
\newtheorem{rem}{Remark}[section]

\numberwithin{equation}{section}
\begin{document}
\begin{center}\large{{\bf{ Weak Henstock-Orlicz space and inclusion properties}}} \\
\vspace{0.5cm}
 
   Hemanta Kalita$^{a}$\footnote{Corresponding author}  and Bipan Hazarika$^{b}$
     
\vspace{0.5cm}
\footnotesize $^a$Department of Mathematics, Assam Don Bosco University, Sonapur 782402, Assam, India.\\
\footnotesize $^{b}$Department of Mathematics, Gauhati University, Guwahati 781014, Assam, India\\
\footnotesize Email:    $^a$hemanta30kalita@gmail.com; $^b$bh\_rgu@yahoo.co.in; bh\_gu@gauhati.ac.in  
\end{center}
\title{}
\author{}
\thanks{\today} 
\begin{abstract}
In this paper we discuss the structure of Henstock-Orlicz space with locally Henstock integrable functions. The weak Henstock-Orlicz spaces on $\mathbb{R}^n$  and some basic properties of the weak Henstock-Orlicz spaces are studied.  We obtain some necessary and sufficient conditions for the inclusion properties of these spaces. 
{\footnotesize {\bf{AMS subject classification \textrm{(2010)}: 26A39,  28B05, 46B03, 46B20, 46B25, 46G10, 46G12}}}
\end{abstract}
\maketitle

\maketitle

\pagestyle{myheadings}
\markboth{\rightline {\scriptsize   Kalita, Hazarika}}
        {\leftline{\scriptsize Weak Henstock-Orlicz spaces}}

\maketitle

    \section{Introduction }
   In 1912 Arnaud Denjoy presented a powerful integral which was able to integrate all finite derivatives and recover their preimitive functions which was not possible in the case of Lebesgue integral. The Denjoy and Perron intergals are generalizations of the Lebesgue integrals that recover a continuous function from its derivative (see \cite{Ra}).  
    J. Kurzweil introduced a generalized version of the Riemann integral (see \cite{Kurzweil}). In 1960's, Henstock made the first systematic study of this new integral. 
     Four years later, while unware of the work of Kurzweil,  Henstock published a paper on integration theory in which he discussed the same integral as Kurzweil. Throughout a series of papers in the sixties Henstock developed a substantial amount of properties of this integral. The definition of this integral as defined by Kurzweil  \cite{Kurzweil} and Henstock  \cite{Henstock} is quite elegant as it is highly reminiscent of the Riemann integral and since a substantial amount of its properties can be developed using Riemann sums and basic epsilon-delta proofs. For the honours of these mathematicians, now a days this integral is called Henstock-Kurzweil integral, also see \cite{Rh}. During late nineties a lot of integration theorist have been studied, the Henstock-Kurzweil integral extensively and consequently the theory of this integral had been highly refined. It should be
pointed out that this integral does not have a standard name at that time. It
is also referred to as the Henstock-Kurzweil integral (in short Henstock integral), the generalized Riemann
integral, and the gauge integral. Since the integrals discussed so far (Riemann,
Lebesgue, Denjoy, Perron) are named for a single person and since Henstock
launched the study of this integral, we are content to call it the Henstock integral (one can see \cite{CD, Piaza} for related works of Henstock-Kurzweil, McShane and Pettis integrals). 
 The  Henstock–Kurzweil integration on Euclidean spaces initiated by Yeong \cite{Lee}. The Orlicz space is the generalization of the $L^p$ space, which was initiated by Z.W. Birnbaum and W. Orlicz. 
 The fundamental properties of Orlicz space with Lebesgue measure found in \cite{Ma}. The theory of Orlicz space which is a more generalized version of $L^p$-space with the help of Young functions and the underlying measure was discussed  in \cite{Mm} (also see \cite{BP}). 
  Nakai used Orlicz spaces in the  application in Harmonic analysis in various ways (see \cite{Nakai, Nakai1}). Several inclusion properties of Orlicz and weak Orlicz spaces are found in \cite{Masta}. Liu Pei De et al. \cite{Lie} discussed about the application of the weak Orlicz spaces in the Harmonic analysis.   By over coming difficulties as $C_0^\infty(\mathbb{R}^n)$  is dense in $L^p (\mathbb{R}^n )$  but not generally dense in $\mathcal{L}^\theta(\mathbb{R}^n)$ Thung in \cite{Tv}, presented  a translation invariant subspace $L^1 (\mathbb{R}^n) \cap \mathcal{L}^\theta (\mathbb{R}^n)$ to be dense in Orlicz space $\mathcal{L}^\theta(\mathbb{R}^ n ) .$  The concept of the Henstock-Orlicz space (in brief $\mathbb{H}-$Orlicz space) was presented by Hazarika and Kalita in \cite{BH} which have some difference thing of the Orlicz spaces such as in $\mathbb{H}-$Orlicz spaces $C_0^\infty$ is dense  but not generally dense in Orlicz spaces. In \cite{HK}, Kalita and Hazarika was investigated the countable additivity of Henstock-dunford integrable functions on $\mathbb{H}-$Orlicz space. The theory of $\mathbb{H}-$Orlicz spaces with vector measure  discussed in the conference paper \cite{HSB}. 
\section*{Objectives of our paper} 
In this paper we  discuss about $\mathbb{H}-$Orlicz space  with a little different settings. In our work we mainly focus on the weak $\mathbb{H}-$Orlicz spaces. Nakai in \cite{Nakai, Nakai1} and Liu Ning et al.  \cite{Liu} defined a class of weak Orlicz function spaces and their basic properties are discussed. The major drawback of the weak Orlicz space is that it is not naturally complete. We motivate to resolve this drawback of the weak Orlicz space. We introduce   $\mathbb{H}-$Orlicz space with locally Henstock integrable functions  to overcome the difficulties of the weak Orlicz space.
   \section{Preliminaries and auxiliary results} In the whole   article, we consider $(\mathbb{R}^n , \Sigma_\infty, \nu_{\infty})$ is an abstract measure space, where   
   $\Sigma_\infty $ is an $\sigma$-algebra of its subsets on which a $\sigma$-additive function $ \nu_{\infty}: \Sigma_\infty \to \mathbb{R}^+ $ is given and $\nu_{\infty}$ 
   is the Lebesgue measure. It is known that a measure space has the finite subset property if for every $A \in \Sigma_\infty$ with $\nu_\infty(A) = \infty$ there exists a family of subsets $\{A_i\}_{i=1}^{\infty} \subset \Sigma_\infty $ with $A_i \subset A;~0 < \nu_\infty(A_i) < \infty$ and $\nu_\infty\left(\bigcup\limits_{i=1}^{\infty}A_i\right) = \infty.$\\
   This give us 
   \[\nu_{\infty}(A)=\left\{\begin{array}{l} 
   0 ,\mbox{ if } A = \emptyset,\\
   +\infty,\mbox{ if } A\neq \emptyset \end{array}\right .\]
   Otherwise it does not restrict the generality of our assumption.
         The space of all Henstock integrable functions defined on $ \mathbb{R}^n$, is denoted by $ HK(\nu_\infty) .$  $ HK(\nu_\infty)$ is a vector space under the usual operations of pointwise addition and scalar multiplication on $\mathbb{R}$ was studied in \cite{Ra,Cs,Bs}.  In the one-dimensional case, Alexiewicz \cite{AL}  has shown that the class  of Henstock  integrable functions, with respect to the  norm 
\[
\left\| h \right\|_{HK}  = \sup_t \left| {\int_{ - \infty }^t {h(s)d(s)} } \right|.
\]
 is a normed space, and it is known that $HK(\mathbb{R})$ is not complete (see \cite{AL}).  
         
  \subsection{\bf Henstock-integral on $\mathbb{R}^n$} 
  The elements of $\mathbb{R}^n$ will be denoted by $z=(z_1,z_2,\dots,z_n).$ An interval in $\mathbb{R}^n$ is a set of the form $\mathcal{J}=[\mathbf{z},\mathbf{w}] := \prod_{i=1}^n[z_i, w_i],$ where
    $-\infty < z_i < w_i <\infty $ for $i = 1,2\dots,n.$ The set $\prod_i^n[z_i, w_i]\subset \mathbb{ R}^n$ is known as a degenerate interval if $z_i =w_i$ for some $i\in \{1,2,\dots,n\}.$ Two intervals $\mathcal{J}=[\mathbf{z},\mathbf{w}], \mathcal{I}= [\mathbf{u}, \mathbf{v}] $ in $\mathbb{R}^n$ are said to be non-overlapping if  $\prod_i^n(z_i, w_i)\cap \prod_i^n(u_i, v_i)$ is empty. The union of
    two intervals in $\mathbb{R}^n$ is  an interval in $\mathbb{R}^n$ (see Lemma 2.1.2 \cite{Lee}).        We know that the space  $\mathbb{R}^n $ equipped with the maximum norm $||.||,$ where 
    $||z|| = \max\limits_{ 1 \leq i \leq n}|z_i|.$  With this norm, we denote the closed ball of $\mathbb{R}^n$ by $B[\mathbf{x},r]=\{\mathbf{x}\in\mathbb{R}^n: ||\mathbf{y}-\mathbf{x}||\le r\},$ whose center is  $\mathbf{x}$ with sides parallel to the co-ordinates axes of length $2r.$ It is a closed interval    for side $i$ about $ x_i$ is in $[z_i,w_i].$ So, let $B[\mathbf{x},r]=[\mathcal{J},\mathbf{x}],$ where $\mathcal{J}= \prod_{i=1}^{n}[z_i,w_i]$, $\mathcal{J} $ is closed interval in $\mathbb{R}^n.$ 
    \begin{defn}
   \cite{BH,Lee} Let $A$ is a compact ball in $\mathbb{R}^n,$ a partition $\mathcal{P}$ of $A$ is a collection $\{(\mathcal{J}_i,x_i): x_i \in \mathcal{J}_i, 1 \leq i \leq m \},$ where $\mathcal{J}_1, \mathcal{J}_2,\dots,\mathcal{J}_m$ are non overlapping closed intervals i.e., $\nu_\infty[\mathcal{J}_i \cap \mathcal{J}_j] = 0, i \neq j$ and $\bigcup\limits_{i=1}^{m} \mathcal{J}_i = A.$
    \end{defn}
     If $\delta $ is a positive function on $A$ we say $\mathcal{P}$ is Henstock $\delta$-partition of $A$ if for each $i$, $\mathcal{J}_i \subset B^{'}(x_i, \delta(x_i)).$ The function $\delta $ is a gauge on $A.$
     \begin{defn}
    \cite{BH,Lee}  A function $h: A \to \mathbb{R}$ is said to be Henstock integrable on $A,$ if there exists a number $L$ such that for any $\varepsilon >0 $ there exist a gauge $\delta$ and Henstock $\delta$-partition on $A$ such that $$\left|\sum\limits_{i=1}^{m}h(x_i)\nu_\infty(\mathcal{J}_i)- L\right| < \varepsilon.$$ 
     \end{defn}
 Now we introduce the concept of locally Henstock-Kurzweil integrable function as follows:
     \begin{defn}\label{bh}
 A measurable function $h:\G \subset \mathbb{R}^n \to \mathbb{R}$ is called locally Henstock-Kurzweil integrable if $h{\chi_\mathcal{K}} \in HK(\mathcal{K})$ for all $\mathcal{K} \subseteq \G$ compact where $\chi_\mathcal{K}$ is the characteristics functions of $\mathcal{K}.$   We denote the set of locally Henstock-Kurzweil integrable functions as $HK_{loc}.$
 \end{defn}
  Recalling that the function $h$  is Henstock integrable on a measurable set $\G \subset \mathbb{R}^n$ if $ h\chi_{\G}$ is Henstock integrable on $\G.$ That is $h \in HK_{loc}(\mathbb{R}^n)$ means $h \in HK(\mathcal{K}),$ where $\mathcal{K}\subseteq \G$ compact.
 Also with easy analogous, $L_{loc}^1(\mathbb{R}^n) \subset HK_{loc}(\mathbb{R}^n).$

              
 \begin{defn}
  \cite{BH, Mm} A function $ \theta: \mathbb{R} \to \mathbb{R}^{+}$ is said to be Young function, so that $\theta(x)= \theta(-x), \theta(0)=0, \theta(x) \to \infty $ as $ x \to \infty,$ but $ \theta(x_0) = +\infty $ for some $ x_0 \in \mathbb{R} $ is permitted. 
  \end{defn}
  	We assume $\mathfrak{F}$ be the class of Young's function $\theta:\mathbb{R}\to \mathbb{R}^+$ is an increasing, bijective, continuous and concave satisfying $\theta(0)=0;~\lim\limits_{t \to 0}\theta(t)=0$ and $\lim\limits_{t \to \infty}\theta(t)=\infty.$ We denote $\theta_1 \Delta \theta_2$ for $\theta_1, \theta_2 \in \mathfrak{F}$ if there is a constant $C>0$ such that $\theta_1(t) \leq \theta_2(Ct)$ for all $t \geq 0.$
    \begin{defn}
        \begin{enumerate}
        \item  A Young function $\theta$ is said to satisfy $\Delta^{\prime}$ if $\lim\limits_{k \to 0}\sup\limits_{t>0}\frac{\theta(kt)}{\theta(t)}=0.$
         \item  An $N$-function $\theta $ is said to satisfy $\Delta_2$-condition if there is a $ k > 0 $ such that $ \theta(2x) \leq k \theta(x)  $ for large values of $ x.$
          \end{enumerate}  
                  \end{defn} 
              If $\theta$ is a convex function on $[0, \infty),$ then $\theta \in \Delta^{\prime}.$ In this article, we do not generally assume that $\theta$ is convex, except we mention it especially.
             \begin{defn}
             \cite[Definition 2.1]{Nakai} Let $\theta$ be a convex function. The Orlicz space is  defined as 
             \begin{eqnarray}
             \mathcal{L}^\theta(\mathbb{R}^n)=\{h \in L_{loc}^{1}(\mathbb{R}^n):~||h||_{L^\theta(\mathbb{R}^n)} < +\infty\},
             \end{eqnarray}
             where $$||h||_{L^\theta(\mathbb{R}^n)}= \inf\left\{ \alpha>0:~(L)\int_{\mathbb{R}^n}\theta\left(\dfrac{|h(x)|}{\alpha}\right)d\nu_\infty(x)\leq 1\right\}~{\textit{~for~some~}} \alpha>0.$$
             \end{defn}                     
   Also one can see \cite{Ma,Nakai1,Mm,BP} for detailed on Orlicz space and it's applications. 
    We recalling few preliminaries of $\mathbb{H}-$Orlicz spaces from \cite{BH}.
 \begin{defn}\cite{BH} 
Let $ (\mathbb{R}^n, \Sigma_\infty, \nu_{\infty})$ be an arbitrary measure space.  Then the space $ \mathcal{H}^{\theta}(\nu_{\infty}) $ of all measurable  functions $ h: \mathbb{R}^n \to \mathbb{R} $ is called $\mathbb{H}-$Orlicz space, which defined as: $$ \mathcal{H}^{\theta}(\nu_{\infty}) = \left\{ h: \mathbb{R}^n \to \mathbb{R} \mbox{~measurable}:  \int_{\mathbb{R}^n}\theta(\alpha h ) d\nu_{\infty} \in HK(\nu_\infty)  \mbox{~for~some~}  \alpha > 0\right\}.$$ 
\end{defn}
 
  \begin{defn}\label{00} \cite{BH}
    The Luxemburg  norm on $\mathcal{H^{\theta}(\nu_{\infty})}$ as follows: $$ {\mathbb{H}_{\theta}} (h)= \inf \left\{ \alpha >0 : HK\int_{\mathbb{R}^n} \theta\left( \frac{h}{\alpha}\right)d\nu_{\infty} \leq 1  \right\}.$$ 
  \end{defn}
  It is understood that $\inf(\emptyset)= + \infty,$  and  $$V= \left\{ h~\mbox{~measurable}: HK\int_{\mcB \subset \mathbb{R}^n}\theta(h)d \nu_\infty \leq 1 \right\},$$ is the gauge of the set.
  \begin{rem}
  The $\mathbb{H}-$Orlicz spaces is Banach spaces with the Luxemberg norm defined as the definition(\ref{00}) as well as Symmetric. This spaces is separable if $\theta \in \Delta_2$ (see \cite{HK}). The separability of $\mathbb{H}-$Orlicz spaces with out $\Delta_2$ conditions of the Young's functions will be an aspect to work in our coming days. 
  \end{rem}
  \section{$\mathbb{H}-$Orlicz spaces with locally Henstock integrable functions}
 Here we discuss about $\mathbb{H}-$Orlicz spaces with a little different settings.
    Before start the concept of $\mathbb{H}-$Orlicz spaces we define the $\mathbb{H}-$Orlicz class state as:
  \begin{defn}\label{de21}
   Let $ \mathcal{H^{-\theta}(\nu_\infty)} $ be the set of all $ h: \mathcal{K} \subset \mathbb{R}^n \to \mathbb{R} $   bounded   measurable with compact support  for $\Sigma_\infty \subset \mathbb{R}^n$  such that $\int_{\mathcal{K}} \theta(|f|)d\nu_{\infty} $ is  Henstock integrable. \\ 
  i.e. $\mathcal{H^{-\theta}(\nu_\infty)} = \{ h  \mbox{~is~bounded~measurable~with~compact~support}: \int_{\mathcal{K}} \theta(|h|)d\nu_{\infty} \in HK(\nu_\infty)\}.$
  \end{defn}
   \begin{thm}
  The space $\mathcal{H^{-\theta}(\nu_\infty)} $ is a linear space if and only if $\theta$ satisfies $\Delta_2$-condition.
  \end{thm}
  \begin{proof}
  The proof is similar as \cite[Theorem 2.2]{BH}.
  \end{proof}
   \begin{prop}\label{prop320}
  	For each $h \in \mathcal{H}^{\theta}(\nu_\infty),$ there is an $\alpha>0$ such that $$ \mathbb{B}_{\theta}= \left\{\alpha h = m \in \mathcal{H}^{-\theta}(\nu_\infty) : HK \int_{\mathcal{K}}\theta(m) d \nu_\infty \leq 1 \right\}$$ is a circled solid subset of $\mathcal{H}^{- \theta}(\nu_\infty).$
  \end{prop}
  \begin{proof} 
  Let $h,m \in \mathcal{H}^{\theta}(\nu_\infty).$ Then there exist $\alpha_0, \beta_0 >0 $ such that $\alpha_0 h, \beta_0 m \in \mathcal{H}^{-\theta}(\nu_\infty).$  Let $\mu_0 = \min(\alpha_0, \beta_0)$. Then for $\mu_0>0 $ and using the known fact of convexity and monotonicity of $\theta,$ we get $$ HK \int_{\mathcal{K}}\theta \left( \frac{\mu_0}{2}(h+m)\right) d \nu_\infty \leq \frac{1}{2}\left[ HK\int_{\mathcal{K}} \theta(\alpha_0 h) d \nu_\infty + HK\int_{\mathcal{K}} \theta(\beta_0 m) d \nu_\infty\right] .$$ Clearly the  right side is Henstock integrable. Since $\frac{\mu_0}{2}> 0 ,$  this gives us $h + m \in \mathcal{H}^{\theta}(\nu_\infty).$  Particularly, with each $ h $ in $\mathcal{H}^{\theta}(\nu_\infty),~2h \in \mathcal{H}^\theta(\nu_\infty) $  and then $n h \in \mathcal{H}^\theta(\nu_\infty)$ for all integers $n >1$, so that  $\gamma_0 h \in \mathcal{H}^\theta(\nu_\infty)$ for any scalar $\gamma_0.$ Therefore the given set is solid and circled. To hold $\gamma_0 h \in \mathcal{H}^\theta(\nu_\infty)$ for some $\gamma_0>0. $ Let $a_n \to 0 $ be arbitrary and set $\gamma_{0_n}= \min(\gamma_0, a_n).$ Then $\gamma_{0_n} \to 0 $ and $\theta(\gamma_{0_n} h) \leq \theta(\gamma_0 h)$ and $\theta(\gamma_{0_n}h) \to 0 $ as $\theta $ is a continuous Young function. Now Dominated Convergence Theorem, give us $HK\int_{\mathcal{K}} \theta(\gamma_{0_n} h) \to 0 $ so that for some $n_0,$ we have $HK\int_{\mathcal{K}}\theta(\gamma_{n_o} h ) d \nu_\infty \leq 1.$ Thus $\gamma_{n_0} h \in \mathbb{B}_{\theta}.$ 
  \end{proof}
  To construct $\mathcal{H}^\theta(\mathbb{R}^n)$ in our setting recall the  known facts that $HK(\mathbb{R}^n) \subset HK_{loc}(\mathbb{R}^n).$   For a Young function $\theta,$ we can state our $\mathbb{H}-$Orlicz space as below:   
  \begin{eqnarray*}
  \mathcal{H}^\theta(\mathbb{R}^n)= \{h \in HK_{loc}: ||f||_{\mathbb{H}_\theta}< +\infty\},
  \end{eqnarray*}
  where $||h||_{\mathbb{H}_\theta(\mathbb{R}^n)}$ is defined  as follows: $$ ||h||_{\mathbb{H}_\theta}= \inf \left\{ \alpha >0 : HK\int_{\mathcal{K}} \theta\left( \frac{h}{\alpha}\right)d\nu_{\infty} \leq 1  \right\},$$  where  $\mathcal{K}$ is a compact subset of $\mathbb{R}^n .$ It is clear that this space is a normed space with respect to the norm $\mathbb{H}_{\theta}$. 

    \begin{thm}\label{th312}
     For each $ h \in \mathcal{H}^{\theta}(\mathbb{R}^n),~m \in \mathcal{H}^{\Phi}(\mathbb{R}^n), $ if the complementary function of $ \theta $ is  $\phi,$ then $$HK\int_{\mathbb{R}^n}|hm|d\nu_{\infty} \leq ||h||_{\mathbb{H}_{\theta}}||m||_{\mathbb{H}_{\phi}} .$$
     \end{thm}
    \begin{thm}
    The classical Orlicz space is a dense subspace of $\mathbb{H}-$Orlicz space as continuous dense embeddings. That is, $\mathcal{L}^\theta(\mathbb{R}^n) \hookrightarrow  \mathcal{H}^\theta(\mathbb{R}^n)$ is continuous dense embeddings.
    \end{thm}
    \begin{proof}
    Let $h \in \mathcal{L}^\theta(\mathbb{R}^n).$ Then $h \in L_{loc}$ with $||f||_{L^\theta(\mathbb{R}^n)}< \infty.$ Then for some $\alpha>0,$ and  a compact $\mathcal{K}\subset \mathbb{R}^n$  we have
    \begin{align*}
  \inf\left\{ HK\int_{\mathcal{K}} \theta\left(\dfrac{h(x)}{\alpha}\right)d\nu_\infty \right\} &\leq \inf\left\{(L)\int_{\mathbb{R}^n}\theta\left(\dfrac{|h(x)|}{\alpha}\right)d\nu_\infty \right\}\\& \leq 1.
   \end{align*}
   So, for some $\alpha>0,~\inf\left\{\alpha>0:~(L)\int_{\mathbb{R}^n} \theta\left(\frac{|h(x)|}{\alpha}\right)d\nu_\infty \leq 1\right\},$ we get the following $$\inf\left\{\alpha>0:~HK\int_{\mathcal{K}} \theta\left(\dfrac{h(x)}{\alpha}\right)d\nu_\infty \leq 1\right\}.$$  Hence $h \in \mathcal{H}^\theta(\mathbb{R}^n)$ with $||h||_{\mathcal{H}^\theta(\mathbb{R}^n)} \leq ||h||_{\mathcal{L}^\theta(\R^n)}.$ Hence the proof.
    \end{proof}
 \begin{thm}
	Suppose $\nu_\infty(\mathcal{K})<\infty$ and $\nu_\infty$ is bounded, then $\mathcal{H}^\theta(\mathbb{R}^n) \hookrightarrow  L^1(\mathbb{R}^n)$ is continuous.
\end{thm}
\begin{proof}
	If $r>0$ and $s>0$ such that for all $p \geq 0,$ we have $\theta(p) \geq rp-s.$ This means $rp \leq \theta(p) +s.$ Let $h \in \mathcal{H}^\theta(\mathbb{R}^n).$ Then for $\alpha>0$ as possible small, we have
	\begin{align*}
	\alpha \int_{\mathcal{K}}|h|d \nu_\infty &\leq \frac{1}{r}\int_{\mathcal{K}}[\theta(\alpha h) +s]d \nu_\infty \\&= \frac{1}{r}\int_{\mathcal{K}} \theta(\alpha h) d\nu_\infty + s \frac{\nu_\infty(\mathcal{K})}{r}\\&< \infty
	\end{align*}
	Therefore $h \in L^1(\mathbb{R}^n).$ Hence $\mathcal{H}^\theta(\mathbb{R}^n) \subset L^1(\mathbb{R}^n).$ 	The inclusion is continuous follows from the similar technique that we used in the second part of the \cite[Lemma 3.4]{BH}. 
\end{proof}
\begin{corollary}
	$\mathcal{H}^\theta(\mathbb{R}^n) \subset HK_{loc}(\mathbb{R}^n).$
\end{corollary}
    \begin{rem}
    If $h \in \mathcal{L}^\theta(\R^n),$ then all theorems that are  true in $\mathcal{L}^\theta(\R^n)$ are so in  $\mathcal{H}^\theta(\mathbb{R}^n).$
    \end{rem}
    \section{Weak-Henstock-Orlicz spaces}
    In this division of the paper, we discuss  about the weak Henstock-Orlicz spaces (In brief weak-$\mathbb{H}-$Orlicz space) and its basic properties. The classical weak Orlicz space is not naturally Banach spaces. We discuss the weak Henstock-Orlicz space is naturally a Banach space. We consider all functions $h$ are $\nu_\infty$ measurable now onward.   $\chi_{\mathcal{B}}$ be the characteristic functions of a measurable $\mathcal{B} \subset \mathbb{R}^n.$  For the  measurable set $\mathcal{B},$ a measurable function $h$ and $t>0,$ let $$\mu_\infty(\mathcal{B}, h,t)=\{x \in \mathcal{B}:~\nu_\infty(h(x))>t\}.$$ If $\mathcal{B}= \mathbb{R}^n,$ we denote as $\nu_\infty(h,t).$ 
    \begin{defn}
    Let $h$ be $\mu_\infty$ measurable function. Then weak-$\mathbb{H}-$Orlicz space is defined as:
    \begin{eqnarray*}
    \mathcal{H}_{w}^{\theta}(\mathbb{R}^n)=\{h \in HK_{loc}(\mathbb{R}^n):~||h||_{H_w^\theta}< + \infty\},
    \end{eqnarray*}
        where $||h||_{H_{w}^\theta}= \inf\left\{~\exists~ \alpha>0:~\sup\limits_{t>0}\theta(t)\nu_\infty\left(\frac{h}{\alpha}, t\right)\leq 1\right\}.$
         \end{defn}  We consider another $\nu_\infty$-measurable set as 
     \begin{eqnarray*}
     B_w^\theta(\mathbb{R}^n)= \left\{ h \in HK_{loc}(\mathbb{R}^n):~\forall~\alpha>0:~\sup\limits_{t>0}\theta(t)\nu_\infty\left(\frac{h}{\alpha}, t\right)< \infty \right\}.
     \end{eqnarray*}
     It is observe that $B_w^\theta(\mathbb{R}^n) \subset \mathcal{H}_w^\theta(\mathbb{R}^n)$ and $B_w^\theta(\mathbb{R}^n)$ is a linear space. If $h=m$ are $\nu_\infty$-a.e. in $\mathcal{H}_w^\theta(\mathbb{R}^n)$, then we call $h=m$ in $\mathcal{H}_w^\theta(\mathbb{R}^n).$
     We define the modular of  function $h$  as   $$\Upsilon_\theta(h)= \sup\limits_{t>0}\theta(t)\nu_\infty\left(\frac{h}{\alpha}, t\right).$$
     \begin{lemma}\label{2.1} 
     \begin{enumerate}
     \item If $h \leq m~$ are $\nu_\infty$-a.e. and $m \in \mathcal{H}_{w}^{\theta},$ then $h \in \mathcal{H}_{w}^{\theta}.$
     \item $\sup\limits_{t>0}\theta\left(\dfrac{t}{||h||_{H_w^\theta}}\right)\nu_\infty(f,t) \leq 1.$
     \item If $\theta \in \Delta_2$ and $||h||_{H_w^\theta(\mathbb{R}^n)} \leq 1 $ then $\Upsilon_\theta(h) \leq ||h||_{H_w^\theta(\mathbb{R}^n)}.$
          \end{enumerate}
     \end{lemma}
     \begin{proof}
      $(1)$ It is straight forward as $h\leq m$ are $\nu_\infty$-a.e. means $|h|\leq |m|$ are $\nu_\infty$-a.e.\\
     $(2)$ From the definition of $||h||_{H_w^\theta(\mathbb{R}^n)}$, there exists $C_k \searrow ||h||_{H_w^\theta(\mathbb{R}^n)}$ such that $\theta\left(\frac{h}{C_k}\right)\nu_\infty(h,t) \leq 1,~\forall ~t>0.$ When $k \to \infty,~ \sup\limits_{t>0}\theta\left(\frac{t}{||h||_{H_w^\theta}}\right)\nu_\infty(h,t) \leq 1.$
     \end{proof}
     \begin{thm}
      If $\theta \in \Delta_{2},~(\mathcal{H}_{w}^{\theta}, ||.||_{H_w^\theta})$ is a quasi-Banach space.
     \end{thm}
     \begin{proof}
      The proof of the result is similar to \cite[Lemma 1.1 (4)]{Lie}.
     \end{proof}
     \begin{thm}\label{2.3}
     For $h,~m\in \mathcal{H}_w^\theta,$ the following inequality holds $$ ||h+m||_{H_w^\theta} \leq ||h||_{H_w^\theta}+||m||_{H_w^\theta}.$$
     \end{thm}
     
     \begin{lemma}
     If $\theta \in \Delta_2$ then the followings are true
     \begin{enumerate}\label{22}
     \item  $\mathcal{H}_w^\theta(\mathbb{R}^n) = B_w^\theta(\mathbb{R}^n).$
     \item If $h_n \in \mathcal{H}_w^\theta(\mathbb{R}^n)$ then $||h_n -h||_{H_w^\theta(\mathbb{R}^n)} \to 0$ if $\Upsilon_\theta(h_n-h) \to 0.$
     \item If $h_n \in \mathcal{H}_w^\theta(\mathbb{R}^n)$ then $\sup\{||h_n||_{H_w^\theta(\mathbb{R}^n)}\}\leq A$ if and only if $\sup\{\Upsilon_\theta(h_n)\} \leq A,$ where $A>0.$
     \end{enumerate}
     \end{lemma}
     \begin{proof}
     $(1)$ Let $\theta \in \Delta_2$ and $
     h \in \mathcal{H}_w^\theta(\mathbb{R}^n).$ Then there exists $\alpha>0$ such that $$\sup\limits_{t>0}\theta\left(\dfrac{t}{\alpha}\right)\nu_\infty(h,t) \leq 1.$$ Hence $\theta(\frac{t}{a})\leq \theta(\frac{t}{\alpha})$ when $a \geq \alpha.$ As $\theta \in \Delta_2,$ so $\theta(\frac{t}{a}) \leq C \theta(\frac{t}{\alpha})$ when $a< \alpha.$ So, $$\sup\limits_{t>0}\theta\left(\frac{t}{a}\right)\nu_\infty(h,t) \leq 1~\mbox{~for~all~}~a>0.$$ Hence $h \in B_w^\theta(\mathbb{R}^n).$\\ 
    The proof of part  $(2) $ follows from the similar technique as \cite[Corollary 2.1]{Lie} and the proof of part $(3)$ follows from \cite[Theorem 2.1 (4)]{Lie}.
     \end{proof}
     The weak Orlicz space is not naturally complete. The weak Orlicz space is complete if the Young function $\theta \in \Delta_2.$ We observe the weak $\mathbb{H}-$Orlicz is complete if the Young function without $\theta \in \Delta_2.$  We discuss now that  the completeness of  weak $\mathbb{H}-$Orlicz space. 
     \begin{lemma}\label{five}
     Let   $(h_n) \in \mathcal{H}_w^\theta(\mathbb{R}^n).$  Then the followings are true
     \begin{enumerate}
     \item If $|| h_n - h||_{H_w^\theta} \to 0,$ then $h_n \to h$ (convergence in measure).
     \item If $0\leq  \inf(h_n) \to h$ is $\nu_\infty-$a.e.  then $||h_n||_{H_w^\theta} \to ||h||_{H_w^\theta}.$
     \end{enumerate}
     \end{lemma}
     \begin{proof}
     (1) The definition of norm of the weak $\mathbb{H}-$Orlicz space, we have 
     \begin{align*}
     ||h_n -h||_{H_w^\theta(\R^n)}&= \inf \left\{\alpha>0:~\sup\limits_{t>0}\theta\left( \frac{h_n-h}{\alpha}\right)\nu_\infty(h_n-h,t)\right\} \to 0
\end{align*} 
gives, $\theta\left(\frac{h_n-h}{\alpha}\right)\nu_\infty(h_n-h,t) \to 0.$  This implies $\nu_\infty(h_n -h) \to 0.$ Therefore $h_n \to h$ in the measure $\nu_\infty.$ \\ 
$(2)$ Let $\inf \nu_\infty(h_n) = \nu_\infty(h)~\nu_\infty-$a.e. Then clearly $h_n \to h $ in the measure $\nu_\infty.$ This gives for all $\varepsilon>0$ there exists a number $n_0 \in \mathbb{N}$ such that $$\nu_\infty(h_n -h) \to 0~\mbox{~for~} n>n_0.$$ Using the property of Young function $\theta,$ we can write
\begin{align*}
\theta\left(\frac{h_n-h}{\alpha}\right)\nu_\infty(h_n-h,t) \to 0&\\ \Rightarrow \sup \theta\left(\frac{h_n -h}{\alpha}\right)\nu_\infty(h_n -h,t) \to 0.
\end{align*} 
That is, $\inf\left\{\alpha>0:~\sup \theta\left(\frac{h_n-h}{\alpha}\right)\nu_\infty(h_n-h,t)\right\}< \varepsilon.$ As  $\varepsilon$  is an arbitrary so $$||h_n-h||_{H_w^\theta(\mathbb{R}^n)} \to 0.$$
     \end{proof}
     \begin{thm}
     $(\mathcal{H}_w^\theta(\R^n), ||.||_{H_w^\theta}(\R^n))$ is a Banach space.
     \end{thm}
     \begin{proof}
     Let $(h_n) \in \mathcal{H}_w^\theta(\mathbb{R}^n)$ be a Cauchy sequence. Then there exists a natural number $ n_0$ such that $\lim\limits_{n,e \to \infty}||h_n-h_e||_{H_w^\theta(\R^n)}=0$ for all  $n,e \geq n_0.$ We know $\lim\limits_{n,e \to n_0}\nu_\infty(h_n-h_e,t)=0$ for all $t>0,$ implies $h_n \to h$ as $\nu_\infty$ measurable. By Riesz's Theorem, there is a subsequence $h_{n_k} \to h,~\nu_\infty-$a.e and $\nu_\infty(h_{n_k}-h_{n_s}) \to \nu_\infty(h_{n_k} -h)$ as $\nu_\infty-$a.e. Let $n_k \geq n_0$ and $h_{n_s} \to \infty,$ then by \cite[Eq 1.1.15 of Rem 1.1.8]{GR} for $\epsilon>0,$
     \begin{align*}
     \theta(\frac{t}{\epsilon})\nu_\infty(h_{n_k}-h,t)&\leq \lim\limits_{n_{s} \to \infty}\theta(\frac{t}{\epsilon})\nu_\infty(h_{n_k}-h,t)\\&\leq 1.
\end{align*}      Using $(2)$ of the  Lemma \ref{five}, $||h_{n_k}-h||_{H_w^\theta(\mathbb{R}^n)}< \varepsilon$ and $\lim\limits_{n_k \to \infty}||h_{n_k}-h||_{H_w^\theta(\R^n)}=0.$ Now from Theorem \ref{2.3}, we can find the following 
     \begin{eqnarray*}
     ||h||_{H_w^\theta(\mathbb{R}^n)} \leq ||h_{n_k}-h||_{H_w^\theta (\mathbb{R}^n)}+||h_{n_k}||_{H_w^\theta(\mathbb{R}^n)}.
     \end{eqnarray*}
     This implies $||h||_{H_w^\theta(\mathbb{R}^n)}\leq ||h_{n_k}||_{H_w^\theta(\mathbb{R}^n)}.$ Using the first part of the Lemma \ref{2.1}, we can conclude that $h \in \mathcal{H}_w^\theta(\mathbb{R}^n).$ So, $h_n \to h $ in $\mathcal{H}_w^\theta(\mathbb{R}^n).$
     \end{proof}
     \begin{thm}
     $\mathcal{L}_w^\theta(\mathbb{R}^n) \hookrightarrow  \mathcal{H}_w^\theta(\mathbb{R}^n)$ is a continuous dense embeddings.
     \end{thm}
     \begin{proof}
     Let $h \in \mathcal{L}_w^\theta(\mathbb{R}^n).$ Then $h \in L_{loc}^{1}(\mathbb{R}^n)$ with $||h||_{H_w^\theta(\mathbb{R}^n)}< \infty.$ This is very obvious that $h \in HK_{loc}.$ Now we find  $||h||_{H_w^\theta(\mathbb{R}^n)}< \infty.$ From the fact that for some $\alpha>0,$
     \begin{align*}
     ||h||_{H_w^\theta(\mathbb{R}^n)} &= \inf\left\{\alpha>0:~\sup\limits_{t>0}\theta(t) \nu_\infty\left(\dfrac{h}{\alpha},t\right)\leq 1\right\}.
     \end{align*}
     We get $\sup\limits_{t>0}\theta(t) \nu_\infty\left(\frac{h}{\alpha},t\right)\leq 1$ also true when $h \in L_{loc}^{1}(\mathbb{R}^n).$ That is, $||h||_{H_w^\theta(\mathbb{R}^n)}\leq ||h||_{L_w^\theta(\mathbb{R}^n)}.$ Hence $\mathcal{L}_w^\theta(\R^n) \hookrightarrow  \mathcal{H}_w^\theta(\R^n)$ is a continuous dense embeddings.
     \end{proof}
     \begin{rem}
      The $\Delta_2$ condition of the Young function $\theta$ is not necessary to proof the Bounded Convergence theorem, Control convergence theorem, Fatou-type convergence theorem, Levi-type convergence theorem, Vitali-type convergence theorem that are proved in the \cite[Section 3]{Liu}, can also be proved  in the weak $\mathbb{H}-$Orlicz space.
     \end{rem}
     \section{Inclusion property of weak $\mathbb{H}-$Orlicz spaces}
     In this section we discuss inclusion properties of  weak $\mathbb{H}-$Orlicz spaces. Before that we find inclusion relations between $\mathbb{H}-$Orlicz space and weak $\mathbb{H}-$Orlicz space in the following theorem as follows:
     \begin{thm}\label{in1}
     Let $\theta$ be a Young function. Then $\mathcal{H}^\theta(\mathbb{R}^n) \subset \mathcal{H}_w^\theta(\mathbb{R}^n)$ for every $h \in \mathcal{H}^\theta(\mathbb{R}^n)$ with $||h||_{H_w^\theta(\mathbb{R}^n)} \leq ||h||_{\mathbb{H}^{\theta}(\mathbb{R}^n)}.$
     \end{thm}
     \begin{proof}
     Let $h \in \mathcal{H}^\theta(\mathbb{R}^n).$ We need to prove $h \in \mathcal{H}_w^\theta(\mathbb{R}^n).$ Let 
     \begin{eqnarray*}
     \mathfrak{A}_{\theta, w}= \left\{\alpha>0:~\sup\limits_{t>0} \theta(t)\nu_\infty\left(\dfrac{h}{\alpha},t\right)\leq 1\right\}
     \end{eqnarray*}
     and \begin{eqnarray*}
     \mathfrak{B}_{\theta, w}= \left\{\alpha>0:~HK\int_{\mathcal{K}} \theta\left(\dfrac{h}{\alpha}\right)d \mu_\infty \leq 1\right\}.
     \end{eqnarray*}
     Clearly, $||h||_{H_w^\theta(\mathbb{R}^n)}= \inf \mathfrak{A}_{\theta, w}$ and $||h||_{\mathbb{H}^\theta(\mathbb{R}^n)}= \inf \mathfrak{B}_{\theta, w}.$ Now for any $\beta \in  \mathfrak{B}_{\theta, w}$ and $t>0$, we have 
     \begin{align*}
     \theta(t)\nu_\infty\left(\frac{h}{\beta},t\right) &\leq HK\int_{\left\{x \in \mathbb{R}^n:~\nu_\infty\left(\frac{h}{\beta},t\right)\right\}}\theta\left(\dfrac{h}{\beta}\right)d \nu_\infty \\&\leq HK\int_{\mathcal{K}}\theta\left(\dfrac{h}{\beta}\right)d\nu_\infty\\&\leq 1.
     \end{align*}
     As $t>0$ is arbitrary, $\sup \limits_{t>0}\theta(t)\nu_\infty\left(\frac{h}{\beta},t\right) \leq 1$ and $\mathfrak{B}_{\theta,w} \leq \mathfrak{A}_{\theta,w}.$  Hence $f \in \mathcal{H}_w^\theta(\mathbb{R}^n)$ with $||h||_{H_w^\theta(\mathbb{R}^n)} \leq ||h||_{\mathbb{H}^{\theta}(\mathbb{R}^n)}.$
     \end{proof}
     \begin{lemma}\label{19july}
     Let $\theta$ be a Young function, $ a \in \mathbb{R}^n$ and $r>0$ be arbitrary. Then 
     \begin{eqnarray*}
     ||\chi_{B(a,r)}||_{H_w^\theta(\mathbb{R}^n)} = \dfrac{1}{\theta^{-1}\left(\dfrac{1}{\nu_\infty(B(a,r))}\right)},
     \end{eqnarray*}
    where  $\nu_\infty(B(a,r)) $ is  the volume of open ball $B(a, r).$
     \end{lemma}
     \begin{proof}
     Theorem \ref{in1} gives $||h||_{H_w^{\theta}(\mathbb{R}^n)} \leq ||h||_{H^{\theta}(\mathbb{R}^n)}.$
     Now,
     \begin{align*}
     ||\chi_{B(a,r)}||_{H^\theta(\mathbb{R}^n)} &= HK\int_{\mathbb{R}^n} \theta\left(\chi_{B(a,r)}{\theta^{-1}\left(\frac{1}{\nu_\infty(B(a,r))}\right)}\right)d\nu_\infty\\&= HK\int_{B(a,r)}\theta\left(\chi_{B(a,r)}\theta^{-1}\left(\frac{1}{\nu_\infty(B(a,r))}\right)\right)d\nu_\infty\\
     &\leq HK\int_{B(a,r)}\frac{1}{\nu_\infty B(a,r)}d\nu_\infty\\&= \frac{1}{\nu_\infty(B(a,r))}HK\int_{B(a,r)} d\nu_\infty\\&=1.
     \end{align*}
     By the definition of $||.||_{H^{\theta}(\mathbb{R}^n)},$ we get $$ ||\chi_{B(a,r)}||_{H_w^\theta\left(\mathbb{R}^n\right)} \leq  \frac{1}{\theta^{-1}(\frac{1}{\nu_\infty(B(a,r))})}.$$ Now we need to prove $$ ||\chi_{B(a,r)}||_{H_w^\theta(\mathbb{R}^n)} \geq  \frac{1}{\theta^{-1}(\frac{1}{\nu_\infty(B(a,r))})}.$$ If possible, let $ ||\chi_{B(a,r)}||_{H_w^\theta(\mathbb{R}^n)} <  \frac{1}{\theta^{-1}(\frac{1}{\nu_\infty(B(a,r))})}.$ Then by the definition of $||.||_{H_w^\theta(\mathbb{R}^n)},$ we find $ ||\chi_{B(a,r)}||_{H_w^\theta(\mathbb{R}^n)} \leq 1.$ This contradicts our assumption and hence \begin{eqnarray*}
     ||\chi_{B(a,r)}||_{H_w^\theta(\mathbb{R}^n)} = \frac{1}{\theta^{-1}(\frac{1}{\mu_\infty(B(a,r))})}.
     \end{eqnarray*}
   \end{proof}
     \begin{thm}
     Let $\theta,~\phi$ be Young functions then the  statements below  are equivalent:
     \begin{enumerate}
     \item $\theta(t) \leq \phi(Ct)$ for every $t>0.$
     \item $\mathcal{H}_w^\phi(\mathbb{R}^n) \subseteq \mathcal{H}_w^\theta(\mathbb{R}^n).$
     \item For every   $h \in \mathcal{H}_w^\theta(\mathbb{R}^n),$ implies $||h||_{H_w^\theta(\mathbb{R}^n)} \leq C||h||_{H_w^\phi(\mathbb{R}^n)}.$
     \end{enumerate}
     \end{thm}
     \begin{proof}
     The proof is similar to the \cite[Theorem 3.3]{Masta}.
     \end{proof}
  \begin{corollary}
 	Let $\theta,~\phi$ be young functions with $\theta(t) \leq \phi(Ct)$ for every $t>0.$ If $h \in \mathcal{H}_{w}^{\theta}(\mathbb{R}^n) $ then $||h||_{H_w^\phi(\mathbb{R}^n)} \leq ||h||_{H_w^\phi(\mathbb{R}^n)}.$
 \end{corollary}
     Now we state a necessary and sufficient condition for the inclusion properties of weak $\mathbb{H}-$Orlicz
spaces generated by concave function presented in the following theorem.
\begin{thm}
Let $\theta_1, \theta_2 \in \mathfrak{F}.$ The followings are equivalent:
\begin{enumerate}
\item $\theta_1 \Delta \theta_2.$
\item $\mathcal{H}_{w}^{\theta_2}(\mathbb{R}^n) \subseteq \mathcal{H}_{w}^{\theta_1}(\mathbb{R}^n).$
\item There exists a constant $C>0$ such that $||h||_{H_{w}^{\theta_1}(\mathbb{R}^n)} \leq C||h||_{H_{w}^{\theta_2}(\mathbb{R}^n)}$ for every $h \in \mathcal{H}_{w}^{\theta_2}(\mathbb{R}^n).$
\end{enumerate}
\end{thm}
\begin{proof}
  $(1) \Rightarrow(2)$ Let $h \in \mathcal{H}_w^{\theta_2}(\mathbb{R}^n).$ We set   
 \begin{eqnarray*}
 \mathbb{A}_{\theta_1}=\left\{\alpha>0:~\sup\limits_{t>0}\theta_1(t)\nu_\infty\left(\frac{h}{\alpha},t\right)\leq 1\right\}
 \end{eqnarray*}
 and 
 \begin{align*}
 \mathbb{A}_{\theta_2}&=\left\{\alpha>0:~\sup\limits_{t>0}\theta_2(Ct)\nu_\infty\left(\frac{h}{\alpha},t\right)\leq 1\right\}\\&= \left\{\alpha>0:~\sup\limits_{t>0}\theta_2(p)\nu_\infty\left(\frac{Ch}{\alpha},p\right)\leq 1\right\}
 \end{align*}
 for $p=Ct.$ If $t>0$ and $\alpha \in \mathbb{A}_{\theta_2},$ then 
 \begin{align*}
 \theta_1(t)\nu_\infty\left(\frac{h}{\alpha},t\right) &\leq \theta_2(Ct)\nu_\infty\left(\frac{h}{\alpha},t\right)\\&= \theta_2(p)\nu_\infty\left(\frac{Ch}{\alpha},t\right)\\&\leq 1.
 \end{align*}
 So, \begin{align*}
 ||h||_{H_w^{\theta_1}(\mathbb{R}^n)} &= \inf \mathbb{A}_{\theta_1}\\&\leq \inf \mathbb{A}_{\theta_2}\\&= C||h||_{H_w^{\theta_2}(\mathbb{R}^n)}.
 \end{align*}
 Therefore $\mathcal{H}_w^{\theta_2}(\mathbb{R}^n) \subseteq \mathcal{H}_w^{\theta_1}(\mathbb{R}^n).$\\
  $(2) \Rightarrow (3)$ Since $(\mathcal{H}_w^{\theta_2}(\mathbb{R}^n), \mathcal{H}_w^{\theta_1}(\mathbb{R}^n))$ are pair of Banach spaces then by \cite[Lemma 3.3]{Krein}, we get the conclusion.\\
   $(3) \Rightarrow(1)$ From the Lemma \ref{19july}, 
 \begin{align*}
   \frac{1}{\theta_{1}^{-1}\left(\frac{1}{\nu_\infty(B(a,r))}\right)}&= ||\chi_{B(a,r)}||_{H_w^\theta(\mathbb{R}^n)} \\&\leq C||\chi_{B(a,r)}||_{H_w^{\theta_2}(\mathbb{R}^n)} \\&= C \frac{1}{\theta_{2}^{-1}\left(\frac{1}{\nu_\infty(B(a,r))}\right)}.
 \end{align*}
 Therefore for any $a \in \mathbb{R}^n$ and $r>0$ we find $$C\theta_{1}^{-1} \left(\frac{1}{\nu_\infty(B(a,r))}\right)\geq \theta_{2}^{-1}\left( \frac{1}{\nu_\infty(B(a,r))}\right).$$ Using the \cite[Lemma 1.1(4)]{Masta}, we get $\theta_1(\frac{1}{\nu_\infty(B(a,r))}) \leq \theta_2(\frac{C}{\nu_\infty(B(a,r))}).$\\
As   $r>0$ is an arbitraty , assuming $t= \frac{1}{\nu_\infty(B(a,r))},$ we get the conclusion $\theta_1(t) \leq \theta_2(Ct)$ and hence $\theta_1 \Delta \theta_2.$
\end{proof}
     
    \section{Declaration}
{\bf Funding:} Not Applicable, the research is not supported by any funding agency.\\
{\bf Conflict of Interest/Competing interests:} The authors declare that the article is free from conflicts of interest.\\
{\bf Availability of data and material:} The article does not contain any data for
analysis.\\
{\bf Code Availability:} Not Applicable.\\
{\bf Author's Contributions:} All the authors have equal contribution for the preparation of the article.

\end{document}